\documentclass[12pt]{article}
\usepackage{amsmath} 
\usepackage{amssymb}

\newcommand{\RR}{{\rm I\hskip -.25em R}}

\begin{document}

\begin{center}
{\bf \large The Killing tensors on an $n$-dimensional manifold with $SL(n,\RR)$-
structure}
\end{center}

\begin{center}
S. E. Stepanov, I. I. Tsyganok \footnote{emails: s.e.stepanov@mail.ru; i.i.tsyganok@mail.ru}
\end{center}

\begin{center}
{\it Finance University, Moscow, Russia}
\end{center}

{\bf Abstract}

In this paper we solve the problem of finding integrals of equations determining 
the Killing tensors on an $n$-dimensional differentiable manifold $M$ endowed 
with an equiaffine $SL(n,\RR)$-structure and discuss possible applications of obtained 
results in Riemannian geometry. 


{\bf AMS:}
53A15; 53A45
\bigskip

{\bf \S1. Introduction}

{\bf 1.1.}
The "structural point of view" of affine differential geometry was introduced 
by K.~Nomizu in 1982 in a lecture at M\"unster University with the title "What 
is Affine Differential Geometry?" (see [No 1]). In the opinion of K.~Nomizu, 
the geometry of a manifold $M$ endowed with an {\it equiaffine structure} is
called affine differential geometry.

In recent years, there has been a new ware of papers devoted to affine 
differential geometry. Today the number of publications (including monographs) 
on affine differential geometry reached a considerable level. The main part of 
these publications is devoted to geometry of hypersurfaces (see [Sch] and [SSV] 
for the history and references). 

{\bf 1.2.} 
In the present paper we solve the problem of finding integrals of equations 
determining the Killing tensors (see [KSCH] for the definitions, properties and 
applications) on an $n$-dimensional differentiable manifold $M$ endowed with an 
equiaffine structure. The paper is a direct continuation of [St 2]. The same 
notations are used here.  

The first of two present theorems proved in our paper is an affine analog of 
the statement published in the paper [St 1], which appeared in the process of 
solving problems in General relativity. 

\bigskip
{\bf \S2. Definitions and results}

{\bf 2.1.} 
In order to clarify the approach to problem of finding integrals of equations 
determining the Killing tensors on an $n$-dimensional differentiable manifold 
$M$ we shall start with a brief introduction to the subject which emphasizes 
the notion of an equiaffine $SL(n,\RR)$-structure.

Let $M$ be a connected differentiable manifold of dimension $n$ $(n>2)$, and 
let $L(M)$ be the corresponding bundle of linear frames with structural group 
$GL(n,\RR)$. We define $SL(n,\RR)$-structure on $M$ as a principal $SL(n,\RR)$-
subbundle of $L(M)$. It is well known that an $SL(n,\RR)$-structure is simply a 
volume element on $M$, i.e. an $n$-form $\eta$ that does not vanishing anywhere 
(see [Ko], Chapter I, \S 2).

We recall the famous problem of the existence of a uniquely determined linear 
connection $\nabla$ reducible to $G$ for each given $G$-structure on $M$ (see 
[C], p. 213). For example, if $M$ is a manifold with a pseudo-Riemannian metric 
$g$ of an arbitrary index $k$, then the bundle $L(M)$ admits a unique linear 
connection $\nabla$ without torsion that is reducible to $O(m,k)$-structure. 
Such a connection is called the {\it Levi-Civita connection}. It is 
characterized by the following condition $\nabla g=0$. 

A linear connection $\nabla$ having zero torsion and reducible to $SL(n,\RR)$ 
is said to be {\it equiaffine} and can be characterized by the following 
equivalent conditions (see [Sch], p 99; [SSV], pp. 57-58):

(1) $\nabla \eta =0$; 

(2) the Ricci tensor $Ric$ of $\nabla$ is symmetric; that means 
$Ric(X,Y)=Ric(Y,X)$ for any vector fields $X,Y\in C^\infty TM$.  

An {\it equiaffine $SL(n,\RR)$-structure} or an {\it equiaffine structure} on 
an $n$-dimensional differentiable manifold $M$ is a pair $(\eta,\nabla)$, where 
$\nabla$ is a linear connection with zero torsion and $\eta$ is a volume 
element which is parallel relative to $\nabla$ (see [No 2], p. 43). 

The curvature tensor $R$ of an equiaffine connection $\nabla$ admits a 
point-wise $SL(n,\RR)$-invariant decomposition of the form 
$R=(n-1)^{-1}[id_M\otimes Ric-Ric\otimes id_M]+W$  where $W$ is the {\it Weyl 
projective curvature tensor} (see [SSV], p. 73-74; [E], \S 40). Then two classes 
of equiaffine structures can be distinguished in accordance with this 
decomposition: the {\it Ricci-flat} equiaffine $SL(n,\RR)$-structures for which 
$Ric=0$, and the {\it equiprojective} $SL(n,\RR)$-structures for which 
$$
R=(n-1)^{-1}[id_M\otimes Ric-Ric\otimes id_M]. \eqno (2.1)
$$

{\bf Remark 1.} 
Recall that a linear connection $\nabla$ with zero torsion is called 
{\it Ricci-flat} if the Ricci tensor $Ric=0$ (see [M]). On the anther hand, a 
connection $\nabla$ is called {\it equiprojective} if the Weyl projective 
curvature tensor $W=0$ (see [Sch], \S 18). In the literature equiprojective 
connections sometimes are called {\it projectively flat} (see, for example, 
[SSV], p. 73). 

An autodiffeomorphism of the manifold $M$ is an automorphism of $SL(n,\RR)$-
structure if and only if it preserves the volume element $\eta$. Let $X$ be a 
vector field on $M$. The function $div X$ defined by the formula 
$(div X)\eta=L_X \eta$ where $L_X$ is the Lie differentiation in the direction 
of the vector field $X$ is called the divergence of $X$ with respect to the $n$-
form $\eta$ (see [KN], Appendix no. 6). Obviously, $X$ is an infinitesimal 
automorphism of an $SL(n,\RR)$-structure if and only if $div X=0$. Such a vector 
field $X$ is said to be {\it solenoidal}. 

For an arbitrary vector field $X$ on $M$ with a linear connection $\nabla$ we 
can introduce the tensor field $A_X=L_X-\nabla_X$ regarded as a field of linear 
endomorphisms of the tangent bundle $TM$. If $M$ is an $n$-dimensional with an 
equiaffine $SL(n,\RR)$-structure then the formula $trace A_X=-div X$ can be 
verified directly (see [KN], Appendix no. 6).  

We have the $SL(n,\RR)$-invariant decomposition 
$A_X=-n^{-1}(div X)id_M+\dot A_X$ at every point $x\in M$. 

Two classes of vector fields on $M$ endowed with an equiaffine $SL(n,\RR)$-structure 
can be distinguished in accordance with this decomposition: the 
{\it solenoidal vector fields} and the {\it concircular vector fields} for 
which, by definition (see [Sc], p. 322; [M]), we have $A_X=-n^{-1}(div X)id_M$.

The integrability conditions of the structure equation $A_X=-n^{-1}(div X)id_M$ 
of the concircular vector field $X$ is the Ricci's identity 
$Y(div X)Z-Z(div X)Y=nR(Y,Z)X$ for any vector fields $Y,Z\in C^\infty TM$ (see 
[E], \S 11). This identity are equivalent to the condition $W(Y,Z)X=0$ for any 
vector fields $Y,Z\in C^\infty TM$. It follows that an equiaffine $SL(n,\RR)$-
structure on an $n$-dimensional manifold $M$ is equiprojective if and only if 
there exist $n$ linearly independent concircular vector fields 
$X_1, X_2, \dots ,X_p$ on $M$ (see also [ST]). This statement is an affine 
analog of the well known fact for the Riemannian manifold $M$ of constant 
sectional curvature (see [F]).

{\bf Remark 2.}  
A pseudo-Riemannian manifold $(M,g)$ with a projectively flat Levi-Civita 
connection $\nabla$ is a manifold of constant section curvature (see [Sch], 
\S 18). Therefore a manifold $M$ endowed with an equiprojective $SL(n,\RR)$-structure 
is an affine analog of a pseudo-Riemannian manifold of constant 
sectional curvature.

{\bf 2.2.} 
We consider an $n$-dimensional manifold $M$ with an equiaffine $SL(n,\RR)$-
structure and denote by $\Lambda^pM$ $(1 \le p \le n-1)$ the $p^{th}$ exterior 
power $\Lambda^p(T^*M)$ of the cotangent bundle $T^*M$ of $M$. Hence 
$C^\infty \Lambda^p M$, the space of all $C^\infty$-sections of $\Lambda^p M$, 
is the space of skew-symmetric covariant tensor fields of degree $p$ 
$(1\le p \le n-1)$.

Let $\gamma : J \subset \RR \to M$ be an arbitrary geodesic on $M$ with affine 
parameter $t\in J$. In this case, we have 
$\nabla_{\frac{d\gamma}{dt}}\frac{d\gamma}{dt}=0$ for the tangent vector 
$\frac{d\gamma}{dt}$ of $\gamma$.

{\bf Definition 1} (see [St 2]). 
{\it A skew-symmetric tensor field $\omega \in C^\infty \Lambda ^p M$ 
$(1\le p \le n-1)$ on an $n$-dimensional manifold $M$ with an equiaffine 
$SL(n,\RR)$-structure is called Killing-Yano tensor of degree $p$ if the tensor 
$i_{\frac{d\gamma}{dt}}\omega :=trace \left(\frac{d\gamma}{dt}\otimes 
\omega\right)$ is parallel along an arbitrary geodesic $\gamma$ on $M$.}

From this definition we conclude that 
$\left(\nabla_{\frac{d\gamma}{dt}}\omega\right)
\left(\frac{d\gamma}{dt},X_2,\dots ,X_p\right)=0$ for any vector fields 
$X_2, \dots , X_p \in C^\infty TM$. Since the geodesic $\gamma$ may be chosen 
arbitrary, the above relation is possible if and only if 
$\nabla \omega \in C^\infty \Lambda ^{p+1}M$, which is equivalent to 
$d\omega=(n+1)\nabla \omega$ for the exterior differential operator 
$d: C^\infty \Lambda ^p M \to C^\infty \Lambda ^{p+1} M$.

Obviously, the set of Killing-Yano tensors of degree $p$ $(1\le p \le n-1)$ 
constitutes an $\RR$-module of tensor fields on $M$, denoted by 
${\bf K}^p (M,\RR)$.

Let $X_1, \dots , X_p$ be $p$ linearly independent concircular vector fields on 
$M$ $(1\le p \le n-1)$. Then direct inspection shows that the tensor field 
$\omega$ of degree $n-p$ dual to the tensor field 
$alt\{X_1\otimes\dots\otimes X_p\}$ relative to the $n$-form $\eta$ is a 
Killing-Yano tensor (see also [St 2]). Therefore on any $n$-manifold $M$ with 
equiprojective $SL(n,\RR)$-structure, there exist at least $n![p!(n-p)!]^{-1}$ 
linearly independent Killing-Yano tensors (see [St 2]). Moreover the following 
theorem is true. 

{\bf Theorem 1.} 
{\it On an $n$-dimensional manifold $M$ endowed with an equiprojective 
$SL(n,\RR)$-structure $(\eta,\nabla)$, there exist a local coordinate system 
$x^1, \dots , x^n$ in which an arbitrary Killing-Yano tensor $\omega$ of degree 
$p$ $(1\le p \le n-1)$ has the components
$$
\omega_{i_1\dots i_p}=e^{(p+1)\psi} (A_{i_0 i_1 \dots i_p}x^{i_0} +
B_{i_1 \dots i_p}) \eqno (2.1)
$$
where $A_{i_0 i_1 \dots i_p}$ and $B_{i_1 \dots i_p}$ are arbitrary constants 
skew-symmetric w.r.t. all their indices and $\psi=(n+1)^{-1}ln(\eta)$.}

From the theorem we conclude that the maximum of linearly independent the 
Killing-Yano tensors is by calculating the number $K^p_n$ of independent 
$A_{i_0 i_1 \dots i_p}$ and $B_{i_1 \dots i_p}$ which exist after accounting 
for the symmetries on the indices. It follows that 
$K^p_n=\frac{(n+1)!}{(p+1)!(n-p)!}$ is the maximum number linearly independent 
the Killing-Yano tensors. 

{\bf Corollary 1.} 
{\it Let $M$ be an $n$-dimensional manifold endowed with an equiprojective 
$SL(n,\RR)$-structure then
$$
dim K^p(M,\RR)=\frac{(n+1)!}{(p+1)!(n-p)!}.
$$}

On our fixed manifold $M$ with an equiaffine $SL(n,\RR)$-structure, we denote 
by $S^pM$ the bundle of symmetric covariant tensor fields of degree $p$ on $M$. 
Hence $C^\infty S^p M$, the space of all $C^\infty$-sections of $S^p M$, is the 
space symmetric covariant tensor fields of degree $p$.

{\bf Definition 2} (see [St 2]). 
{\it A symmetric tensor field $\varphi\in C^\infty S^p M$ on an $n$-dimensional 
manifold $M$ with an equiaffine $SL(n,\RR)$-structure is called Killing tensor 
of degree $p$ if 
$\varphi\left(\frac{d\gamma}{dt},\dots ,\frac{d\gamma}{dt}\right)=const$ along 
an arbitrary geodesic $\gamma$ on $M$.} 

Let $\varphi\left(\frac{d\gamma}{dt},\dots ,\frac{d\gamma}{dt}\right)=const$ 
along an arbitrary geodesic $\gamma$ on $M$ and hence $\varphi$ is a Killing 
tensor. Then the above relation is possible if and only if 
$\delta^*\varphi:=\sum_{cicl}\{\nabla \varphi\}=0$ where for the local 
components $\nabla_{i_0}\varphi_{i_1 \dots i_p}$ of $\nabla \varphi$ we define 
the sum $\sum_{cicl}\{\nabla_{i_0}\varphi_{i_1 \dots i_p}\}$ as the sum of the 
terms obtained by a cyclic permutation of the indices $i_0, i_1, \dots,i_p$.

Obviously, the set of Killing tensors constitutes an $\RR$-module of tensor 
fields on $M$, denoted by ${\bf T}^p(M,\RR)$.

Let $M$ be an $n$-dimensional manifold endowed with an equiprojective 
$SL(n,\RR)$-structure $(\eta,\nabla)$, and $\omega_1, \dots , \omega_p$ be $p$ 
linearly independent Killing-Yano tensors of degree $1$ on $M$. Then direct 
inspection shows that the tensor field 
$\varphi :=sym\{\omega_1\otimes \dots \otimes \omega_p\}$ is a Killing tensor 
of degree $p$. Therefore on any $n$-manifold $M$ with equiprojective $SL(n,\RR)$-
structure, there exist at least $(n+p-1)![p!(n-1)!]^{-1}$ linearly independent 
Killing tensors (see also [SS]). Moreover the following theorem is true. 

{\bf Theorem 2.}
{\it On an $n$-dimensional manifold $M$ endowed with an equiprojective 
$SL(n,\RR)$-structure $(\eta,\nabla)$, there exist a local coordinate system 
$x^1, \dots , x^n$ in which the components $\varphi_{i_1 \dots i_p}$ of 
an arbitrary Killing tensor $\varphi$ of degree $p$ can be expressed in the 
form of an $p^{th}$ degree polynomial in the $x^i$'s
$$
\varphi_{i_1\dots i_p}=e^{2p\psi}
\sum\limits^p_{q=0}A_{i_1\dots i_p j_1\dots j_q}x^{j_1}\dots x^{j_q} \eqno (2.2)
$$
where the coefficients $A_{i_1\dots i_p j_1\dots j_q}$ are constant and 
symmetric in the set of indices $i_1,\dots ,i_p$ and the set of indices 
$j_1,\dots ,j_q$. In addition to these properties the coefficients 
$A_{i_1\dots i_p j_1\dots j_q}$ have the following symmetries 
$$
\sum\nolimits_{cicl}\{A_{i_1\dots i_p j_1\dots j_{p-s}}\}_{j_{p-s+1}}=0 \eqno (2.3)
$$
for $s=1,\dots , p-1$ and  
$$
\sum_{cicl}\{A_{i_1\dots i_p j_1}\}=0. \eqno (2.4)
$$ }                                               

From the theorem we conclude that the maximum number of linearly independent 
the Killing tensors is obtained by calculating the number $T^p_n$ of independent 
$A_{i_1\dots i_p j_1\dots j_q}$ $(q=0,1,\dots , n)$ which exist after 
accounting for the symmetries on the indices the dependence relations (2.3) and 
(2.4). By [KL] it follows that 
$$
T^p_n=\frac{p(p+1)^2(p+2)^2 \dots (m+p-1)^2(m+p)}{(p+1)!p!}
$$
is the maximum number linearly independent the Killing-Yano tensors. Then we 
have the following proposition.

{\bf Corollary 2.} 
{\it Let $M$ be an $n$-dimensional manifold endowed with an equiprojective  
$SL(n,\RR)$-structure then
$$
dim T^p(M,\RR)=\frac{p(p+1)^2(p+2)^2 \dots (m+p-1)^2(m+p)}{p!(p+1)!}.
$$}

\bigskip
{\bf \S3. Proofs of theorems}

{\bf 3.1.} 
We let $f:\bar M\to M$ denote the mapping of an $\bar n$-dimensional manifold 
$\bar M$ endowed with an equiaffine $SL(\bar n,\RR)$-structure onto another an 
$n$-dimensional manifold $M$ endowed with an equiaffine $SL(n,\RR)$-structure, 
and let $f_*$ be the differential of this mapping. For any covariant tensor 
field $\omega$ on $M$, we can then define the covariant tensor field 
$f^* \omega$ on $\bar M$, where $f^*$ is the transformation transposed to the 
transformation $f_*$.

If $dim \bar M=dim M=n$ and $f:\bar M\to M$ is a projective diffeomorphism, i.e., 
a mapping that transforms an arbitrary geodesic in $\bar M$ into a geodesic in 
$M$, then we have the following lemma.

{\bf Lemma 1.}  
{\it Let $f:\bar M\to M$ be a projective diffeomorphism of $n$-dimensional 
manifolds endowed with the equiaffine $SL(n,\RR)$-structures 
$(\bar \eta, \bar \nabla)$ and $(\eta,\nabla)$ respectively. Then for an 
arbitrary Killing-Yano tensor $\omega$ of degree $p$ $(1\le p \le n-1)$ on the 
manifold $M$ the tensor field $\bar \omega = e^{-(p+1)\psi}(f^*\omega)$ with 
$\psi=(n+1)^{-1}ln(\eta/\bar \eta)$ will be the Killing-Yano tensor of degree 
$p$ on the manifold $\bar M$.}

{\bf Proof.} 
It is known that the diffeomorphism $f:\bar M\to M$ can be realized following 
the principle of equality of the local coordinates 
$\bar x^1=x^1, \dots ,\bar x^n=x^n$ at the corresponding points $\bar x$ and 
$x=f(\bar x)$ of these manifolds. In this case, we have the equalities (see 
[Sch], \S 18; [M] and [YI])
$$
\Gamma^k_{ij}=\bar \Gamma^k_{ij}+\psi_i\delta^k_j+\psi_j\delta^k_i \eqno(3.1)
$$
for the objects $\Gamma^k_{ij}$ and $\bar \Gamma^k_{ij}$ of the a equiaffine 
connections $\nabla$ and $\bar \nabla$ in the coordinate system 
$x^1, \dots , x^n$ that is common w. r. t. the mapping $f:\bar M\to M$, and 
for the gradient \linebreak $\psi_j=(n+1)^{-1}\partial_j ln[\eta/\bar \eta]$.

Equalities (3.1) imply that the mapping $f^{-1}$, which in inverse to the 
projective diffeomorphism $f:\bar M\to M$, is a projective mapping.

We set $\omega_{i_1 \dots i_p}$ be the local components of a Killing-Yano 
tensor $\omega$ of degree $p$ $(1\le p\le n-1)$ arbitrary defined on the 
manifold $M$; by definition, these components satisfy the equations
$$
\nabla_{i_0}\omega_{i_1\dots i_p}+\nabla_{i_1}\omega_{i_0\dots i_p}=0. 
\eqno(3.2)
$$
From equalities (3.2) we find directly that the components
$$
\bar \omega_{i_1 \dots i_p}=e^{-(p+1)\psi}\omega_{i_1 \dots i_p} \eqno (3.3)
$$
of the tensor field $\bar \omega=e^{-(p+1)\psi}(f^*\omega)$ satisfy the 
equations
$$
\bar\nabla_{i_0}\bar\omega_{i_1\dots i_p}+
\bar\nabla_{i_1}\bar\omega_{i_0\dots i_p}=0 . \eqno (3.4)
$$
Hence, the tensor field $\bar \omega$ is a Killing-Yano tensor of degree $p$ 
$(1\le p \le n-1)$ on the manifold $\bar M$.

{\bf 3.2.} 
Let ${\bf A}^n$ be an $n$-dimensional affine space with a volume element given 
by the determinant: $det (e_1, \dots , e_n)=1$, where $\{e_1, \dots , e_n\}$ is 
the standard basis of the underlying vector space for ${\bf A}^n$. We denote by 
$\nabla$ the standard linear connection in ${\bf A}^n$ relative to which the 
volume element "det" is parallel (see [No 2]; [SSV], p. 10). 

Let $f:\bar M\to {\bf A}^n$  be a projective diffeomorphism from a manifold 
$\bar M$ endowed with equiaffine $SL(n,\RR)$-structure onto an affine space 
${\bf A}^n$ endowed with standard equiaffine S$L(n,\RR)$-structure. It is well 
known that manifolds endowed with equiprojective $SL(n,\RR)$-structures and 
only these manifolds are projectively diffeomorphic to an affine space 
${\bf A}^n$ (see [Sch], \S 18; [M]) therefore in our case the $SL(n,\RR)$-
structure of the manifold $\bar M$ must be an equiprojective structure.

If ${\bf A}^n$ is an affine space with the Cartesian system of coordinates 
$\bar x_1, \dots , \bar x^n$ then the components $\bar \omega_{i_1\dots i_p}$ of 
the Killing-Yano tensor $\bar \omega$ of degree $p$ $((1\le p\le n-1))$ in 
equation (3.4) must now satisfy
$$
\partial_j \bar \omega_{i i_1\dots i_p}+
\partial_i \bar \omega_{j i_1\dots i_p}=0  \eqno (3.5)
$$
where $\partial_j=\frac{\partial}{\partial x^j}$. From (3.5) we conclude the 
following equations
$$
\partial_k \partial_j \bar \omega_{ii_1\dots i_p}+
\partial_k \partial_i \bar \omega_{ji_1\dots i_p}=0; \eqno (3.6)
$$
$$
\partial_j \partial_i \bar \omega_{ki_1\dots i_p}+
\partial_j \partial_k \bar \omega_{ii_1\dots i_p}=0; \eqno (3.7)
$$
$$
\partial_i \partial_k \bar \omega_{ji_1\dots i_p}+
\partial_i \partial_j \bar \omega_{ki_1\dots i_p}=0. \eqno (3.8)
$$
From (3.6), (3.7), (3.8) we find
$$
\partial_k \partial_j \bar \omega_{i_1i_2\dots i_p}=0, \eqno (3.9)
$$
by using identities $\frac{\partial ^2 h}{\partial \bar x^k\partial \bar x^j}=
\frac{\partial ^2 h}{\partial \bar x^j \partial \bar x^k}$ which are carried 
out for an arbitrary smooth function $h:{\bf A}^n\to \RR$. The integrals of 
equations (3.9) take the form
$$
\bar \omega_{i_1\dots i_p}=A_{i_0i_1\dots i_p}\bar x^{i_0}+
B_{i_1\dots i_p} \eqno (3.10) 
$$
for any skew-symmetric constants $A_{i_0i_1\dots i_p}$ and $B_{i_1\dots i_p}$ 
(see also [SS] and [St 3]). Taking the components (3.10) of the Killing-Yano 
tensor $\bar \omega$ in ${\bf A}^n$ and using Lemma 1, we can formulate 
Theorem 1.

{\bf 3.3.} 
Let $\bar M$ be a manifold of dimension $n$ endowed with the equiaffine 
$SL(n,\RR)$-structure $(\bar \eta, \bar \nabla)$ and $M$ be a manifold of some 
dimension endowed with the equiaffine $SL(n,\RR)$-structure $(\eta,\nabla)$. 
Let there is given a projective diffeomorphism $f:\bar M\to M$, then we have 
the following lemma.

{\bf Lemma 2.} 
{\it Let $f:\bar M\to M$ be a projective diffeomorphism of $n$-dimensional 
manifolds endowed with the equiaffine $SL(n,\RR)$-structures 
$(\bar \eta, \bar \nabla)$ and $(\eta,\nabla)$ respectively. Then for an 
arbitrary Killing tensor $\varphi$ of degree $p$ on the manifold $M$ the tensor 
field $\bar \varphi=e^{-2p\psi}(f^*\varphi)$ with 
$\psi=(n+1)^{-1}ln(\eta/\bar \eta)$ will be the Killing tensor of degree $p$ on 
the manifold $\bar M$.}

{\bf Proof.} 
We set $\varphi_{i_1\dots i_p}$ to be components of the Killing tensor $\varphi$ 
arbitrary defined on the manifold $M$; by definition, these components satisfy 
the following equations $\sum_{cicl}\{\nabla_{i_0}\varphi_{i_1\dots i_p}\}=0$. 
Then we find directly that the components $\bar \varphi_{i_1\dots i_p}=
e^{-2p\psi}\varphi_{i_1\dots i_p}$ of the tensor $\bar \varphi=
e^{-2p\psi}\varphi$ satisfy the equations
$$
\sum\nolimits_{cicl}\{\bar\nabla_{i_0}\bar\varphi_{i_1\dots i_p}\}=
e^{-2p\psi}\sum\nolimits_{cicl}\{\nabla_{i_0}\varphi_{i_1\dots i_p}\}=0. \eqno (3.11)
$$
From (3.11) we conclude that the tensor field $\bar \varphi$ is a Killing tensor 
of degree $p$ on the manifold $\bar M$.

{\bf 3.4.} 
It follows from Nijenhuis (see [Ni]) that in an $n$-dimensional affine space 
${\bf A}^n$ the components $\bar \varphi_{i_1\dots i_p}$ of the Killing tensor 
$\bar \varphi$ of degree $p$ can be expressed in the form of an $p^{th}$ degree 
polynomial in the $\bar x^i$'s 
$$
\varphi_{i_1\dots i_p}=e^{-2p\psi}\sum^p_{q=0}
A_{i_1\dots i_pj_1\dots j_q}\bar x^{j_1}\dots \bar x^{j_q}. \eqno (3.12)
$$
The coefficients $A_{i_1\dots i_pj_1\dots j_q}$ are constant and symmetric in 
the set of indices $i_1,\dots ,i_p$ and the set of indices $j_1,\dots ,j_q$. 
In addition to these properties the coefficients $A_{i_1\dots i_pj_1\dots j_q}$ 
have the following symmetries 
$\sum_{cicl}\{A_{i_1\dots i_pj_1\dots j_{p-s}}\}_{j_{p-s+1}}=0$ for $s=1,\dots, p-1$ 
and $\sum_{cicl}\{A_{i_1\dots i_pj_1}=0\}$. Taking the components (3.12) of the 
Killing tensor $\bar \varphi$ in ${\bf A}^n$ and using Lemma 2, we can 
formulate Theorem 2.

\bigskip
{\bf \S4.  Applications to Riemannian geometry}

{\bf 4.1.} 
Let $(M,g)$ be a pseudo-Riemannian manifold of dimensional $n$. Then from the 
present theorems $1$ and $2$ we conclude that an arbitrary Killing vector 
$\omega$ has the following local covariant components 
$\omega_i=e^{2\psi}(A_{ik}x^k+B_i)$ where $\psi=[2(n+1)]^{-1}ln|det g|$, 
$A$'s and $B$'s are constants and $A_{ik}+A_{ki}=0$ (see also [St 1]). It 
follows that the group of infinitesimal isometric transformations has 
$\frac{1}{2}n(n+1)$ parameters (see also [E], \S 71).

{\bf 4.2.} 
Following [T] and [Ka], a skew-symmetric covariant tensor field $\vartheta$ of 
degree $p$ $(1\le p\le n-1)$ is called a conformal Killing tensor if 
$\vartheta \in ker D$ for
$$
D=\nabla - \frac{1}{p+1}d-\frac{1}{n-p+1}g\wedge d^*
$$ 
where $d^*$ is the codifferential operator 
$d^*:C^\infty \Lambda ^{p+1} M\to C^\infty \Lambda ^p M$ and 
$$
(g\wedge d^* \vartheta)_{i_0i_1\dots i_p}=
\sum_{a=1}^{p}(-1)^{a+1}g_{i_0i_a}(d^*\vartheta)_{i_1\dots {\hat i_a}\dots i_p}.
$$
Obviously, the set of conformal Killing tensors constitutes an vector space of 
tensor fields on $(M,g)$, denoted by ${\bf C}^p(M,\RR)$ (see [St 5]). If a 
conformal Killing tensor $\vartheta$ belongs to $ker d^*$, then it is a 
Killing-Yano tensor. On the other hand, if a conformal Killing tensor 
$\vartheta$ belongs to $ker d$, it is called a closed conformal Killing tensor 
or a planar tensor (see [St 4] and [St 5]). We denote the vector space of these 
tensors by ${\bf P}^p(M,\RR)$.

By [Ka] on an arbitrary $n$-dimensional pseudo-Riemannian manifold $(M,g)$ of constant 
nonzero sectional curvature $C$ $(C\ne 0)$ the vector space ${\bf C}^p(M,\RR)$ of 
conformal Killing tensors is decomposed uniquely in the form 
$$
{\bf C}^p(M,\RR)={\bf K}^p(M,\RR)\oplus {\bf P}^p(M,\RR). \eqno (4.1)
$$
From (4.1) we conclude that any conformal Killing tensor $\vartheta$ of degree 
$p$ is decomposed uniquely in the form $\vartheta=\omega+\theta$ where $\omega$ 
is a Killing-Yano tensor of degree $p$ and $\theta$ is a closed conformal 
Killing tensor of degree $p$. 

Following theorem $1$, on an $n$-dimensional pseudo-Riemannian manifold $(M,g)$ 
of constant nonzero sectional curvature $C$ ($C\ne 0$) there is a local 
coordinate system $x^1,\dots , x^n$ in which an arbitrary Killing-Yano tensor 
$\omega$ of degree $p$ $(2\le p \le n-1)$ has the components 
$$
\omega_{i_1\dots i_p}=
e^{(p+1)\psi}(A_{i_0i_1\dots i_p}x^{i_0}+B_{i_1\dots i_p}) \eqno (4.2)
$$
where $\psi=[2(n+1)]^{-1}ln|det g|$, 
$\psi_k=\frac{\partial \psi}{\partial x^k}$ and $A_{i_0 i_1\dots i_p}$, 
$B_{i_1\dots i_p}$ are arbitrary skew-symmetric constants. On the other hand, 
by [St 3] on a pseudo-Riemannian manifold $(M,g)$ of constant nonzero curvature 
$C$ $(C\ne 0)$ the components $\theta_{i_1\dots i_p}$ of a closed conformal 
Killing tensor $\theta$ of degree $p$ $(1\le p\le n-1)$ can be found from the 
equations 
$$
\theta_{i_1 i_2 \dots i_p}=-\frac{1}{pC}\nabla_{i_1}\omega_{i_2\dots i_p} \eqno (4.3)
$$
where $\nabla_{i_1}\omega_{i_2\dots i_p}=\partial_{i_1}\omega_{i_2\dots i_p}-
\omega_{k\dots i_p}\Gamma^k_{i_2 i_1}-\dots - \omega_{i_2\dots k}\Gamma^k_{i_pi_1}$ 
is the expression for the covariant derivative $\nabla\omega$ of the 
Killing-Yano tensor of degree $p-1$. Moreover, by virtue of (3.1) on a 
pseudo-Riemannian manifold $(M,g)$ of constant curvature $C$ $(C\ne 0)$ the 
Christoffel symbols $\Gamma^k_{ij}$ have the following form 
$\Gamma^k_{ij}=\psi_i\delta^k_j+\psi_j\delta^k_i$ (see also [St 1]). Therefore, 
we can deduce from (4.2) and (4.3) that 
$$
\theta_{i_1\dots i_p}=-\frac{1}{C}e^{p\psi}(\psi_{[i_1}A_{|k|i_2\dots i_p]}x^k+
\psi_{[i_1}B_{i_2\dots i_p]}+\frac{1}{p}A_{i_1i_2\dots i_p}).
$$
Consequently we have 

{\bf Theorem 3.} 
{\it On an $n$-dimensional pseudo-Riemannian manifold $(M,g)$ of constant 
nonzero sectional curvature $C$ $(C\ne 0)$ there is a local coordinate system 
$x^1,\dots ,x^n$ in which an arbitrary conformal Killing tensor $\vartheta$ of 
degree $p$ $(2\le p\le n-1)$ has the components
$$
\vartheta_{i_1\dots i_p}=e^{(p+1)\psi}(A_{ki_1\dots i_p}x^k+B_{i-1\dots i_p})-
\frac{1}{C}e^{p\psi}(\psi_{[i_1}C_{|k|i_2\dots i_p]}x^k+
\psi_{[i_1}D_{i_2\dots i_p]}+\frac{1}{p}C_{i_1 i_2\dots i_p})
$$
where $\psi=[2(n+1)]^{-1}ln|det g|$, $\psi_k=\frac{\partial \psi}{\partial x^k}$ 
and $A_{i_0 i_1\dots i_p}$, $B_{i_1\dots i_p}$, $C_{i_1\dots i_p}$ and 
$D_{i_1\dots i_p}$ are arbitrary skew-symmetric constants.}

{\bf Remark 3.} For a conformal Killing vector field, see K. Yano and T. 
Nagano [YN].

\bigskip
{\bf References}

[C] Chern S.S. The geometry of G-structures. {\it Bull. Amer. Math. Soc.} {\bf 72}(1966), 167-219.

[E] Eisenhart L.P. {\it Riemannian geometry}, Princeton Univ. Press, Princeton, NJ, 1949.

[F] Fulton C.M. Parallel vector fields. {\it Proc. Amer. Math. Soc.} {\bf 16}(1965), 136-137.

[KL] Katzin G.H., Levine J. Note on the number of linearly independent $m^{th}$-order first integrals space of constant curvature. {\it Tensor}, {\bf 19}(1968), no. 1, 42-44.

[Ka] Kashiwada T. On konformal Killing tensor. {\it Natural Science Report, Ochanomizu University} {\bf 19}(1968), no. 2, 67-74.

[Ko] Kobayashi Sh. {\it Transformation groups in differential geometry.} Erbeb. Math.     Grenzgeb, Bd. 70, Springer-Verlag, New York-Heidelberg, 1972.

[KN] Kobayashi and Nomizu K. {\it Foundations of differential geometry.} Vol. 1, Intersience, New York-London, 1963.

[KSCH] Kramer D., Stephani H., Mac Callum M.A.H. and Herit E. {\it Exact solutions of Einstein's field equations.} Cambridge Univ. Press, Cambridge, 1980.

[M] Mikes J. Geodesic mapping of affine-connected and Riemannian spaces. {\it J. Math. Sci.}, New York {\bf 78}(1996), no. 3, 311-333.

[Ni] Nijenhuis A. A note on first integrals of geodesics. {\it Proc. Kon. Ned. Acad. Van. Wetens.,} Ser. A {\bf 52}(1967), 141-145.

[No 1] Nomizu K. {\it What is affine differential geometry?} Differential Geometry Meeting Univ. Munster 1982. Tagungsbericht, 1982, 42-43.

[No 2] Nomizu K. On completeness in affine differential geometry. {\it Geometriae Dedicata}, {\bf 20}(1986), no. 1, 43-49.

[Sc] Schouten J.A. {\it Ricci-calculus}. Grundlehren Math. Wiss., Bd. {\bf 10}, Springer-Verlag, Berlin etc., 1954.

[Sch] P.A. and A.P. Schirokow. {\it Affine Differentialgeometrie.} Teubner, Leipzig, 1962.

[SSV] Simon U., Schwenk-Schellschmidt A., Viesel H. {\it Introduction to the affine differential geometry of hypersurfaces.} - Science University of Tokyo, Japan, 1991.

[St 1] Stepanov S.E. The Killing-Yano tensor. {\it Theoretical and Mathematical Physics} {\bf 134}(2003), no. 3, 333-338.

[St 2] Stepanov S.E. The Bochner technique for an $m$-dimensional compact manifold with $SL(m,\RR)$-structure. {\it St. Petersburg Mathematical Journal} {\bf 10}(1999), no. 4, 703-714.

[St 3] Stepanov S.E. On conformal Killing 2-form of the electromagnetic field. {\it Journal of  Geometry and Physics} {\bf 33}(2000), 191-209.

[St 4] Stepanov S.E. A class of closed forms and special Maxwell's equations. {\it Tensor, N.S.} {\bf 58}(1997), 233-242.

[St 5] Stepanov S.E. The vector space of conformal Killing forms on a Riemannian manifold. {\it Journal of  Mathematical  Sciences.} {\bf 110}(2002), no. 4, 2892-2906.

[SS] Stepanov S.E., Smol'nikova M.V. Fundamental differential operators of orders one on exterior and symmetric forms. {\it Russ. Math. J.} {\bf 46}(2002), no.11, 51-56.

[ST] Stepanov S.E., Tsyganok I.I. {\it Vector fields in a manifolds with equiaffine connections.} Webs and Quasigroups, Tver Univ. Press, Tver, 1993, 70-77.

[T] Tachibana S.-I. On conformal Killing tensor in a Riemannian space. {\it Tohoku Math. Journal} {\bf 21}(199), 56-64.

[YI] Yano K., Ishihara Sh. Harmonic and relative affine mappings. {\it J. Differential Geometry} {\bf 10}(1975), 501-509.

[YN] Yano K., Nagano T. Some theorems on projective and conformal transformations. {\it Ind. Math.} {\bf 14}(1957), 451-458.

\end{document}